\documentclass[prl,twocolumn,twoside,superscriptaddress,showpacs,floatfix]{revtex4}
\usepackage{amsmath,graphicx,multirow,natbib,hyperref}
\hypersetup{colorlinks,citecolor=black,filecolor=black,linkcolor=black,urlcolor=black}

\begin{document}

\title{Kuramoto model with uniformly spaced frequencies: \\ Finite-$N$ asymptotics of the locking threshold}
\author{Bertrand Ottino-L\"{o}ffler and Steven H. Strogatz}
\affiliation{Center for Applied Mathematics, Cornell University, Ithaca, New York 14853}
\date{\today}



\begin{abstract}
We study phase locking in the Kuramoto model of coupled oscillators in the special case where the number of oscillators, $N$, is large but finite, and the oscillators' natural frequencies are evenly spaced on a given interval. In this case, stable phase-locked solutions are known to exist if and only if the frequency interval is narrower than a certain critical width, called the locking threshold. For infinite $N$, the exact value of the locking threshold was calculated 30 years ago; however, the leading corrections to it for finite $N$ have remained unsolved analytically. Here we derive an asymptotic formula for the locking threshold when $N \gg 1$. The leading correction to the infinite-$N$ result scales like either $N^{-3/2}$ or $N^{-1}$, depending on whether the frequencies are evenly spaced according to a midpoint rule or an endpoint rule. These scaling laws agree with numerical results obtained by Paz\'{o} [Phys. Rev. E 72, 046211 (2005)]. Moreover, our analysis yields the exact prefactors in the scaling laws, which also match the numerics.  
\end{abstract}

\maketitle


\section{Introduction}

In 1975, Kuramoto proposed an elegant model of coupled nonlinear oscillators, now known as the Kuramoto model \cite{kuramoto75, kuramoto84}. Since then the model has been applied to a wide range of physical, biological, chemical, social, and technological systems, and its analysis has stimulated theoretical work in nonlinear dynamics, statistical physics, network science, control theory, and pure mathematics. For reviews, see Refs.~\cite{strogatz00, pikovsky03, strogatz03, acebron05, dorfler14, pikovsky15, rodrigues16}. 

The governing equations for the Kuramoto model are 
\begin{alignat}{1}\label{BaseKura}
\dot\theta_i = \omega_i + \frac{K}{N} \sum_{j =1}^{N} \sin(\theta_j - \theta_i), 
\end{alignat}
for $i = 1, \ldots, N$, where $\theta_i$ is the phase of oscillator $i$ and $\omega_i$ is its natural frequency. Inspired by Winfree's work on self-synchronizing systems of biological oscillators \cite{winfree67}, Kuramoto restricted attention to attractive coupling, $K>0$, and assumed that the $\omega_j$ were randomly distributed across the population according to a prescribed probability distribution $g(\omega)$, which he took to be unimodal and symmetric about its mean. Without loss of generality the mean frequency can be set to 0 by going into a rotating frame, and the coupling $K$ can be set to $K=1$ by rescaling time. Both of these normalizations will be assumed in what follows.

One adjustable parameter remains: the characteristic width $\gamma$ of the frequency distribution. When $\gamma$ is sufficiently large, numerical simulations show that the oscillators behave incoherently and run at their natural frequencies. At the other extreme, when $\gamma = 0$ the oscillators have identical frequencies and approach a perfectly synchronized solution with $\theta_i(t)=\theta_j(t)$ for all $i, j$ and $t$. Kuramoto's achievement was to analyze the dynamics of the model in between these two extremes. He solved the model in the continuum limit $N \rightarrow \infty$ and obtained a number of beautiful results~\cite{kuramoto75, kuramoto84}, opening up a fruitful line of research for many subsequent studies (reviewed in ~\cite{strogatz00, pikovsky03, strogatz03, acebron05}). In particular, he showed that as $\gamma$ is decreased from large values, a bifurcation takes place at a critical value $\gamma_c$. This bifurcation gives rise to a branch of \emph{partially synchronized} states:  the oscillators with natural frequencies near the mean (= 0) lock to that frequency, forming a synchronized pack, while the others run at different frequencies from the pack and from each other.  Along this branch, the order parameter $r$, defined by
\begin{equation} \label{OrderParam}
re^{i\psi} = \frac{1}{N} \sum_{j=1}^{N} e^{i\theta_j},
\end{equation}
grows continuously from 0 as $\gamma$ decreases through $\gamma_c$.  
Most strikingly, Kuramoto showed that the onset of partial synchronization at $\gamma_c$ is analogous to a second-order phase transition, and he derived exact results for the order parameter along the partially synchronized branch. 

For finite $N$, however, much less is known. The most important advances have come in three areas:  finite-size corrections to the critical coupling at the phase transition, and finite-size scaling laws for the dynamical fluctuations of the order parameter just past the transition \cite{daido87, daido89, daido90, hong05, hong07, hong15};  finite-$N$ corrections to the model's kinetic theory \cite{hildebrand07, buice07}; and analysis of the model's phase-locked states and their bifurcations \cite{ermentrout85, vanhemmen93, jadbabaie04, maistrenko04, aeyels04, mirollo05, pazo05, quinn07, verwoerd08, verwoerd11, dorfler11}.  

Our work in this paper was motivated by an open problem in the third vein, about the threshold for phase locking. To explain what this means and why the question is interesting, we briefly review two definitions and two prior studies. By a \emph{phase-locked state}, we mean a solution of \eqref{BaseKura} that satisfies $\dot\theta_i(t) = \dot\theta_j(t)$ for all $i, j$ and $t$. Equivalently, a phase-locked state is a solution in which all the mutual phase differences $\theta_i(t) - \theta_j(t)$ are unchanging in time.  The earliest results about phase locking in the Kuramoto model were obtained by Ermentrout~\cite{ermentrout85}. In 1985, he analyzed the infinite-$N$ version of \eqref{BaseKura} for frequency distributions $g(\omega)$ supported on a bounded interval $[-\gamma, \gamma]$.  The restriction to distributions without tails was necessary to avoid trivialities; otherwise phase locking would always be impossible in the infinite-$N$ limit. Ermentrout found that phase-locked solutions were possible if and only if $\gamma$ was less than a critical value $\gamma_L$, which we will refer to as the \emph{locking threshold}.  In particular, for a uniform distribution of frequencies on the interval $[-\gamma, \gamma]$, Ermentrout proved that the infinite-$N$ limit of \eqref{BaseKura} has phase-locked states precisely when $\gamma \le \gamma_L^{\infty}=\pi/4$.  Hence the locking threshold is $\pi/4$ in this case. 

In 2005, Paz\'{o}~\cite{pazo05} analyzed a finite-$N$ counterpart of the same problem.  He considered two schemes for arranging $N$ frequencies evenly on the interval $[-\gamma, \gamma]$, thereby providing two different deterministic approximations to the same uniform distribution. One scheme, which we refer to as the \textit{midpoint rule}, is given by 
\begin{equation}
\omega_{j} = \gamma \left(-1 + \frac{2 j-1 }{N} \right), \label{midpoint}
\end{equation}
for $j = 1, \ldots, N$. This rule comes from dividing the interval $[-\gamma, \gamma]$ into $N$ equal subintervals and then placing the frequencies at their midpoints. The second scheme, which we call the \textit{endpoint rule}, is given by 
\begin{equation}
\omega_{j} = \gamma \left(-1 + 2\frac{j-1 }{N-1}\right), \label{endpoint}
\end{equation}
for $j = 1, \ldots, N$.
This rule also spaces $N$ frequencies evenly on $[-\gamma, \gamma]$, but starts at the endpoint $\omega_1=-\gamma$ and continues in equal steps to the other endpoint,  $\omega_N=\gamma$. By computing the locking thresholds $\gamma_L$ numerically for both rules, Paz\'{o}~\cite{pazo05} observed that for $N \gg 1$,
\begin{equation}
\left| \gamma_L - \gamma_L^{\infty} \right| \propto N^{-\mu} \label{midpoint_pazo}
\end{equation}
with $\mu \approx 1.502$ for the midpoint rule \eqref{midpoint}, and 
\begin{equation}
\left| \gamma_L - \gamma_L^{\infty} \right| \propto N^{-1} \label{endpoint_pazo}
\end{equation}
for the endpoint rule \eqref{endpoint}. He then derived \eqref{endpoint_pazo} analytically under the assumption that \eqref{midpoint_pazo} is correct. The derivation used the facts that $\gamma_L^{\infty}=\pi/4$ and
\begin{equation}
 \gamma_L(\text{endpoint}) = \left(1-\frac{1}{N}\right)\gamma_L(\text{midpoint}),  \notag
\end{equation}
which holds because the midpoint rule maps onto the endpoint rule if we rescale $\gamma$ by a factor of $1-1/N$. 
To see that this is the right scaling factor, notice that when $\gamma= \gamma_L$ the maximal frequency satisfies
\begin{equation} 
\omega_N= \gamma_L
\label{gamma_omega_endpoint}
\end{equation}
 for the endpoint rule \eqref{endpoint}, whereas 
 \begin{equation}
 \omega_N = \gamma_L (1-1/N) 
 \label{gamma_omega_midpoint}
 \end{equation}
 for the midpoint rule \eqref{midpoint}. 
 Thus 
\begin{alignat}{1}
\ \gamma_L(\text{endpoint}) &= \left(1-\frac{1}{N}\right)\left( \frac{\pi}{4} +  O(N^{-\mu})\right) \notag \\
&= \frac{\pi}{4} -\frac{\pi}{4} N^{-1} + O(N^{-\mu}). \notag
\end{alignat}
The unsolved problem, however, was to prove \eqref{midpoint_pazo} itself. 

In this paper we derive the leading asymptotic behavior of the locking threshold $\gamma_L$ for the midpoint and endpoint rules. Besides accounting for the scaling exponents of $-3/2$ and $-1$ seen numerically, the asymptotics also yield exact expressions for the prefactors. We prove that
\begin{equation}
\gamma_L = \frac{\pi}{4}+ 4 \zeta\left(-\frac{1}{2} ,\frac{C_1}{2} \right) N^{-3/2} + O(N^{-2}) \notag
\end{equation}
for the midpoint rule \eqref{midpoint} and
\begin{equation}
\gamma_L = \frac{\pi}{4} - \frac{\pi}{4}N^{-1} + 4 \zeta\left(-\frac{1}{2} ,\frac{C_1}{2} \right) N^{-3/2} + O(N^{-2}) \notag
\end{equation}
for the endpoint rule \eqref{endpoint}. Here 
\begin{equation}
\zeta(s, q) = \sum_{n=0}^{\infty} (n+q)^{-s}
\label{Hurwitz}
\end{equation} 
is the Hurwitz zeta function and $C_1 \approx 0.605$ is the QRS constant, defined by Bailey et al.~\cite{bailey09} as the unique zero of $\zeta(1/2, z/2)$ in the interval $0<z<2$.

\section{The locking threshold}

\subsection{Background} 
We recall a few standard results about phase locking in the finite-$N$ Kuramoto model. It it convenient to rewrite the governing equations \eqref{BaseKura} as 
\begin{equation}
\dot\theta_i = \omega_i + Kr\sin(\psi-\theta_i), \notag
\end{equation}
where $r$ and $\psi$ are defined in \eqref{OrderParam}. Then, because we are assuming that the mean $\omega_i$ has been normalized to zero by going into a suitable rotating frame, phase-locked solutions satisfy $\dot\theta_i =0$ for $i=1, \ldots, N$. By rotational symmetry, we can simplify the system further by choosing coordinates such that $\psi=0$. Thus $\omega_i = Kr\sin\theta_i$. Rescaling time so that $K=1$, we obtain 
\begin{equation}
\omega_i= r\sin \theta_i.
\label{omega_vs_theta}
\end{equation}
With these choices, the order parameter simplifies as well; since $\psi=0$,  the real part of \eqref{OrderParam} gives 
\begin{equation}
r=\frac{1}{N} \sum_{j=1}^{N} \cos \theta_j,
\notag
\end{equation}
or equivalently, 
\begin{equation}
r=\langle\cos \theta_j \rangle
\label{self-consistency}
\end{equation}
where the angle brackets denote an average over all oscillators. 

The condition that determines the locking threshold for finite $N$ has been derived by several authors \cite{aeyels04, mirollo05, quinn07, verwoerd08, dorfler11}. It is given implicitly by 
\begin{equation}
2\left \langle\cos \theta_j \right \rangle=\left \langle 1/\cos \theta_j \right \rangle.
\label{saddle-node_condition}
\end{equation}
Equation \eqref{saddle-node_condition} in turn yields a formula for the maximum phase that an oscillator can reach before stable phase locking is lost. As one would expect intuitively, that maximal phase is achieved by the oscillator with the maximal frequency, $\omega_N$. Quinn et al.~\cite{quinn07} showed that corresponding phase $\theta_N$ at the locking threshold is given implicitly by 
\begin{equation}
2\left \langle \sqrt{1-\nu_j^2 \sin^2 \theta_N} \right \rangle = \left \langle  \frac{1} {\sqrt{1-\nu_j^2 \sin^2 \theta_N}}      \right \rangle.
\label{S_N_implicit}
\end{equation}
Here the normalized frequencies $\nu_j$ are defined by 
\begin{equation}
\nu_j = \omega_j/\omega_N
\label{nu_definition}
\end{equation}
and satisfy $-1 \le \nu_j \le 1$; in fact, for both the midpoint rule \eqref{midpoint} and the endpoint rule \eqref{endpoint}, 
\begin{equation}
\nu_j = 1 - 2 \frac{j-1}{N-1}
\label{nu_formula}
\end{equation}
for $j=1, \ldots, N$. By solving \eqref{S_N_implicit} numerically for these $\nu_j$, and assuming $N \gg 1$, Quinn et al.~\cite{quinn07} found that the maximal phase at the locking threshold satisfies
\begin{equation}
\sin \theta_N \approx 1 - \frac{C_1}{N}
\notag
\end{equation}
where $C_1 \approx 0.605443657$.
Bailey et al.~\cite{bailey09} then took the calculations out to 1500 digits and, by a skillful asymptotic analysis, proved that 
\begin{equation}
\sin \theta_N \sim 1 - \frac{C_1}{N} - \frac{C_2}{N^2} + O(N^{-3})
\label{bailey_result}
\end{equation}
where  $C_1$ is the unique zero of the Hurwitz zeta function $\zeta(1/2, z/2)$ in the interval $0<z<2$. Furthermore, they proved that the second coefficient in the series is related to the first by 
\begin{equation}
C_2 = C_1 - C_{1}^2 - 30 \frac{\zeta(-1/2, C_1/2)}{ \zeta(3/2, C_1/2)}. 
\label{miracle}
\end{equation}
Numerically, $C_2 \approx -0.104$. In later work they calculated the coefficients $C_3$ and $C_4$ of the third and fourth terms of the asymptotic series as well~\cite{durgin08}. 

\subsection{Formula for the maximal frequency}

The results above can be leveraged to give an exact formula for the locking threshold $\gamma_L^N$. From \eqref{omega_vs_theta} the maximal phase $\theta_N$ is related to the maximal frequency $\omega_N$ via 
\begin{equation}
\omega_N = r \sin \theta_N.
\label{omega_N}
\end{equation}
 In turn, $\omega_N$ is related to the locking threshold $\gamma_L^N$ via Eqs.~\eqref{gamma_omega_endpoint} and \eqref{gamma_omega_midpoint}. So if we can derive a formula for $r$ in terms of $\sin \theta_N$ (whose asymptotics are known), we will have opened a pathway to understanding the large-$N$ asymptotics of the locking threshold.
 
 To obtain the desired formula for $r$, first combine Eqs.~\eqref{self-consistency} and \eqref{saddle-node_condition} to get
 \begin{equation}
 r=\frac{1}{2} \left \langle \frac{1}{\cos \theta_j} \right \rangle. \notag
 \end{equation}
Then express $\cos\theta_j$ in terms of $\sin \theta_N$, as follows. Divide $\omega_j=r\sin \theta_j$ by $\omega_N=r\sin \theta_N$ to get $\omega_j/\omega_N=\sin \theta_j/\sin \theta_N$. Hence \eqref{nu_definition} implies
\begin{equation}
\sin\theta_j=\nu_j \sin\theta_N
\label{sj}
\end{equation}
and so 
\begin{equation}
r=\frac{1}{2} \left \langle \frac{1}{\sqrt{1-\nu_j^2 \sin^2 \theta_N}} \right \rangle
\end{equation}
where the $\nu_j$ are given by \eqref{nu_formula}. Thus \eqref{omega_N} becomes 
\begin{equation}
\omega_N=\frac{1}{2} \sin\theta_N \left \langle \frac{1}{\sqrt{1-\nu_j^2 \sin^2 \theta_N}} \right \rangle.
\label{omega_N_sum}
\end{equation}
The technical challenge, then, is to analyze the asymptotic behavior of the sum \eqref{omega_N_sum}, where $\sin \theta_N$ is given by \eqref{bailey_result} and $\nu_j = 1 - 2 (j-1)/(N-1)$ with $j=1, \ldots,  N \gg 1$.

\subsection{Rewriting the sum}

Let us recast \eqref{omega_N_sum} as a Riemann sum. From \eqref{sj} we get 
\begin{equation}
\omega_N=\frac{1}{2} \sin\theta_N \left \langle \frac{1}{\sqrt{1- \sin^2 \theta_j}} \right \rangle.
\label{omega_N_sum_with_sj}
\end{equation}
Recall that $\nu_j = \omega_j/\omega_N = \sin \theta_j / \sin \theta_N$ so the spacing between consecutive values of $\nu_j$ is 
\begin{equation}
\Delta \nu = \frac{\Delta s}{\sin \theta_N}, \notag
\end{equation}
where $s$ is shorthand for $\sin \theta$. On the other hand, the explicit formula \eqref{nu_formula} for $\nu_j$ gives $\Delta \nu = 2/(N-1)$.  Equating these two expressions for $\Delta \nu$ gives 
\begin{equation}
\Delta s = \frac{2 \sin \theta_N}{N-1} \notag
\end{equation}
(which makes sense since $s=\sin \theta$ changes from $-\sin\theta_N$ to $\sin\theta_N$ in $N-1$ steps). Replacing the factor $\sin \theta_N$ in \eqref{omega_N_sum_with_sj} with $(N-1) \Delta s/2$ then gives the desired Riemann sum:
\begin{align}
\omega_N &= \frac{1}{2} \left( \frac{N-1}{2} \right) \Delta s \left[ \frac{1}{N} \sum_{j=1}^{N} \frac{1}{\sqrt{1- \sin^2 \theta_j}} \right] \notag\\
                  &=\frac{1}{4} \left( 1-\frac{1}{N} \right) \sum_{j=1}^{N} \frac{\Delta s}{\sqrt{1- s_j^2}} 
\label{Riemann_sum_s}
\end{align}
where 
\begin{equation}
s_j = \sin \theta_j.  \notag
\end{equation}
As $N\to\infty$, the Riemann sum \eqref{Riemann_sum_s} for the maximal frequency $\omega_N$ converges to $(1/4) \int_{-1}^{1} \left(1 - s^2 \right)^{-1/2} ds = \pi/4$, which agrees with Ermentrout's result for the locking threshold in the continuum limit \cite{ermentrout85}. Notice that this integral has integrable singularities at its endpoints $s=\pm 1$. The counterparts of those singularities have be dealt with carefully when estimating the sum for finite $N$.

\section{Estimating the sum}

\subsection{Changing the notation}

Before beginning the asymptotic analysis, it turns out to be helpful to reindex the sum \eqref{Riemann_sum_s} and to revise the notation accordingly. Let $k=j-1$ and $M=N-1$, so that the sum now runs from $k=0$ to $k=M$ instead of $j=1$ to $j = N$.  Replace the dependent variable $s_j$ with its reindexed version, denoted $u_k$ and defined by $u_k \equiv u_{j-1} \equiv s_j$.  Then \eqref{Riemann_sum_s} becomes 
\begin{equation}
\omega_M= \frac{1}{4} \left( \frac{M}{M+1} \right) \alpha_M,  \notag
\end{equation}
where 
\begin{alignat}{1}
\alpha_M = \sum_{k = 0}^{M} \frac{\Delta u} {\sqrt{1 - u_{k}^2}}   \label{alpha_M}.
\end{alignat}
\noindent
Let $A_M(k)$ denote the summands of $\alpha_M$:
\begin{equation}
 A_M(k) = \frac{\Delta u} {\sqrt{1 - u_{k}^2}}.
\label{A_M}
 \end{equation} 
In \eqref{A_M} the mesh points that were formerly $s_j$ have now become 
 \begin{equation}
 u_k = -S_M + k \Delta u
 \label{u_k}
 \end{equation}
 and their spacing is  
 \begin{equation} 
 \Delta u = \Delta s = 2S_M/M, 
 \label{Delta_u}
 \end{equation} 
 where $M=N-1$ and
 $S_M$ is the transformed version of $\sin\theta_N$. From \eqref{bailey_result},
 \begin{equation} 
 S_M  \sim 1 - \frac{C_1}{M+1} - \frac{C_2}{(M+1)^2} + O(M^{-3}).
 \end{equation} 
 Expanding $S_M$ in inverse powers of $M$ yields 
\begin{equation} 
S_M \sim  1 - \frac{C_1}{M} - \frac{C_2-C_1}{M^2} + O(M^{-3}).
 \label{SM}
 \end{equation}

\subsection{Strategy} 

What makes the sum in \eqref{alpha_M} difficult to estimate is that $S_M \rightarrow 1$ and $\Delta u \rightarrow 0$ as $M \rightarrow \infty$; hence more and more of the mesh points $u_k = -S_M + k \Delta u$ approach the integrable singularities at the endpoints $u=\pm1$ as $M \rightarrow \infty$. That is why naive approximation methods for estimating sums based on integrals fail for this problem. A further complication is that the mesh points have a nonstandard dependence on $M$ through $S_M$. 

To make progress, we need to understand and isolate the behavior of a typical summand $A_M(k)$ as we approach the singularities at $u = \pm 1$. In particular, near $u= -1$ we want to find an approximation $D^{-}_{M}(k)$ such that $A_M(k) = D^{-}_{M}(k) + O(M^{-2})$ for small $k$. Likewise, we seek an approximation $D^{+}_{M}(k)$ that captures the dominant behavior of $A_M(k)$ near $u=+1$, where $k=O(M)$. Actually, because our problem has left-right symmetry, it satisfies identities like $u_k^2=u_{M-k}^2$ and $A_M(k)= A_M(M-k)$. Therefore it suffices to study the dominant behavior near one of the endpoints, say $u=-1$, and then double its contribution to any relevant sums by invoking the symmetry.

Our strategy thus involves three steps: (1) Isolate, but do not yet evaluate, the dominant contributions $D^{-}_{M}(k)$ and $D^{+}_{M}(k)$ coming from the \textit{fringes} near the endpoint singularities. (2)  Subtract off the dominant terms and estimate the remaining sum, which we call the \textit{bulk}. This portion of the sum is sufficiently well behaved that it can be approximated by an integral. (3) Only now evaluate the dominant asymptotic behavior in the fringes by summing $D^{-}_{M}(k)$ and $D^{+}_{M}(k)$.  

\subsection{Isolating the dominant terms}

The first step is to substitute \eqref{u_k}, \eqref{Delta_u}, and \eqref{SM} into \eqref{A_M}. Then expand the resulting expression for $A_M(k)$ in powers of $M$ while holding $k$ fixed. The result is
\begin{alignat}{1}
A_M(k) =& M^{-1/2} \left[ \frac{1}{(C_1/2+k)^{1/2}} \right] \notag \\
&+ M^{-3/2} \left[ \frac{-3 C_1^2+2 C_1-2 C_2+4 k^2 } {8 \left(C_1/2+ k\right){}^{3/2}}  \right] \notag \\
&+ O(M^{-5/2}). \notag
\end{alignat}
As we will see below, it turns out to be more helpful to write $A_M(k)$ in the following equivalent form, which has the nice property that the dependence on $k$ occurs solely though the variable $(C_1/2+k)/M$:
\begin{alignat}{1}
A_M(k) =
& \frac{1}{M} \left[1 - \frac{C_1}{2M} \right] \left( \frac{C_1/2+k}{M}\right)^{-1/2} \notag\\
& + \frac{1}{2M} \left(\frac{C_1/2+k}{M} \right)^{1/2} \notag \\
& + \frac{1}{M^3} \left[\frac{C_1 - C_1^2- C_2}{4}\right]  \left(\frac{C_1/2+k}{M} \right)^{-3/2} \notag \\
&+ O(M^{-5/2}). \notag
\end{alignat}
Note that because the error term is smaller than $O(M^{-2})$, we can identify the first three terms above as the desired approximation $D^{-}_{M}(k)$ to the summand $A_M(k)$, for $k$ fixed and $M \gg 1$. Thus
\begin{alignat}{1}
D^{-}_{M}(k) = & \frac{1}{M} \left[1 - \frac{C_1}{2M} \right] \left( \frac{C_1/2+k}{M}\right)^{-1/2} \notag \\
&+ \frac{1}{2M} \left(\frac{C_1/2+k}{M} \right)^{1/2} \notag \\
& + \frac{1}{M^3} \left[\frac{C_1 - C_1^2- C_2}{4}\right]  \left(\frac{C_1/2+k}{M} \right)^{-3/2}.
\label{D_M} 
\end{alignat}

\begin{figure}[!htpb]
        \includegraphics[width=0.9\linewidth]{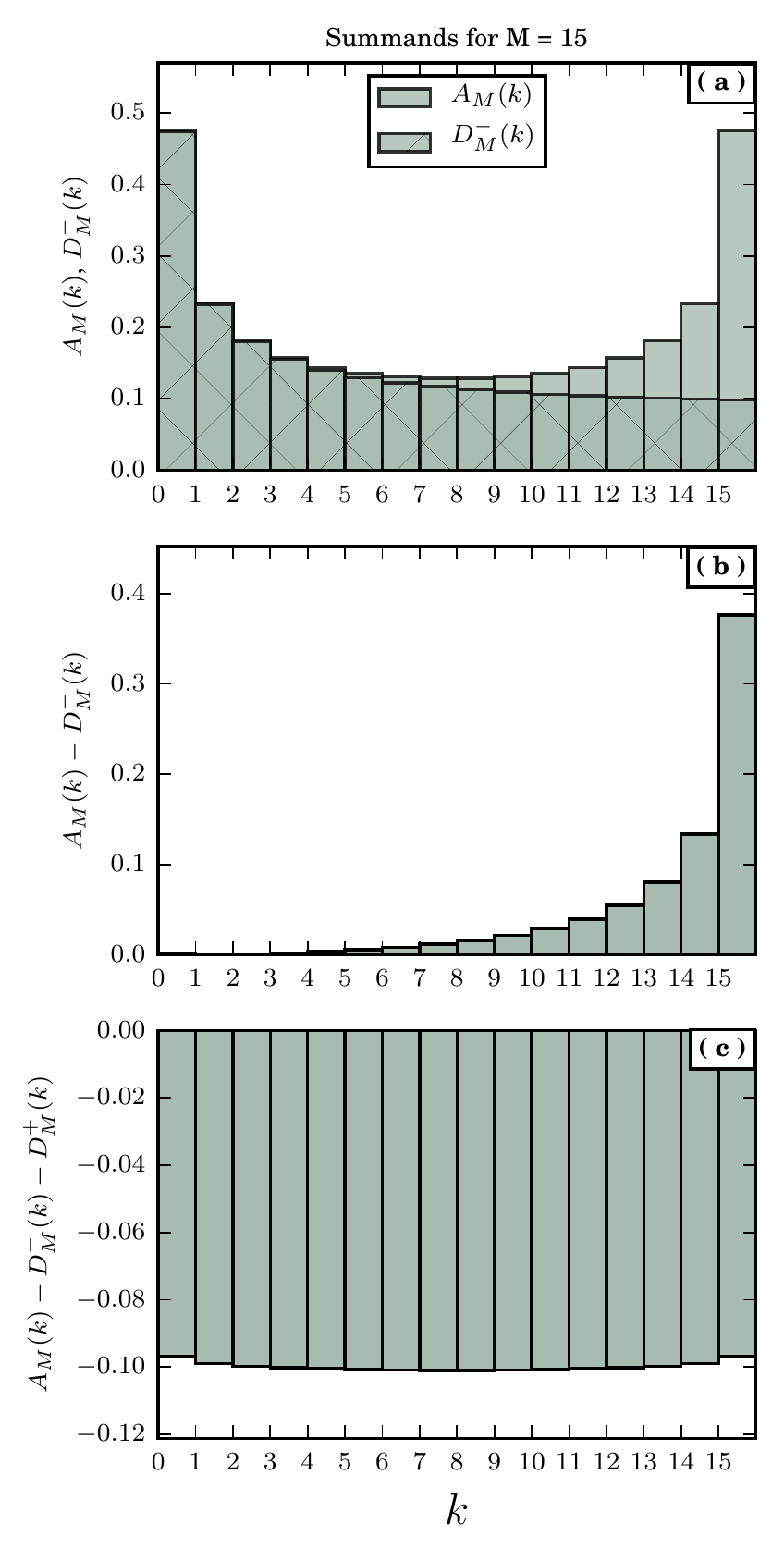}
        \caption{\label{k_pictures} (Color online) Summands $A_M(k)$ and their approximations $D_{M}^{-}(k)$ and $D_{M}^{+}(k)$ near the left and right edges, for $M=15$. (a) Values of $A_M(k)$ (solid rectangles) lie above those of $D_{M}^{-}(k)$ (hatched rectangles). Note that the approximation is most accurate for small values of $k$.  (b) When $D_{M}^{-}(k)$ is subtracted from $A_M(k)$, a near cancellation occurs near the left edge of the region. (c) When both $D_{M}^{-}(k)$ and $D_{M}^{+}(k)$ are subtracted from $A_M(k)$, the remainder is well behaved for all $k$. }
\end{figure}

Now break the sum $\alpha_M$ into two parts, which we denote as $B_M$ ($B$ for bulk) and $F_M$ ($F$ for fringe):
\begin{alignat}{1}
\alpha_M  \notag
=& \sum_{k=0}^{M}  \left[ A_M(k) - D_{M}^{-}(k) - D_{M}^{+}(k)\right]\notag \\
&+ \sum_{k=0}^{M}\left[ D_{M}^{-}(k) + D_{M}^{+}(k) \right] \notag \\
=: &  B_M + F_M. \label{two_part_sum_1} 
\end{alignat}

When the left-right symmetry is applied, this becomes
\begin{alignat}{1}
\alpha_M  \notag
=& \sum_{k=0}^{M} \left[ A_M(k) - 2 D_{M}^{-}(k) \right]+ \sum_{k=0}^{M}\left[ 2 D_{M}^{-}(k)  \right] \notag \\
= &  B_M + F_M. \label{two_part_sum_2} 
\end{alignat}

\subsection{Bulk}
To estimate the bulk sum, $B_M$, we examine its summand $A_M(k) - D_{M}^{-}(k) - D_{M}^{+}(k)$, as shown in Fig.~\ref{k_pictures}.   By construction, $A_M(k) - D_{M}^{-}(k$) decays at least as fast as $O(M^{-2})$  as $k$ approaches 0, and $A_M(k) - D_{M}^{+}(k)$ similarly decays as $k$ approaches $M$.  Moreover, each of $D_{M}^{\pm}(k)$ has a single singularity which cancels that at corresponding endpoint in $A_M(k)$. So in this sense, the sum $B_M$ does not inherit the singular behavior that the summand of $\alpha_M$ once had, meaning that we can safely replace it with a corresponding midpoint integral at the cost of an error no larger than $O(M^{-2})$.  

\begin{figure}[!htpb]
        \includegraphics[width=0.9\linewidth]{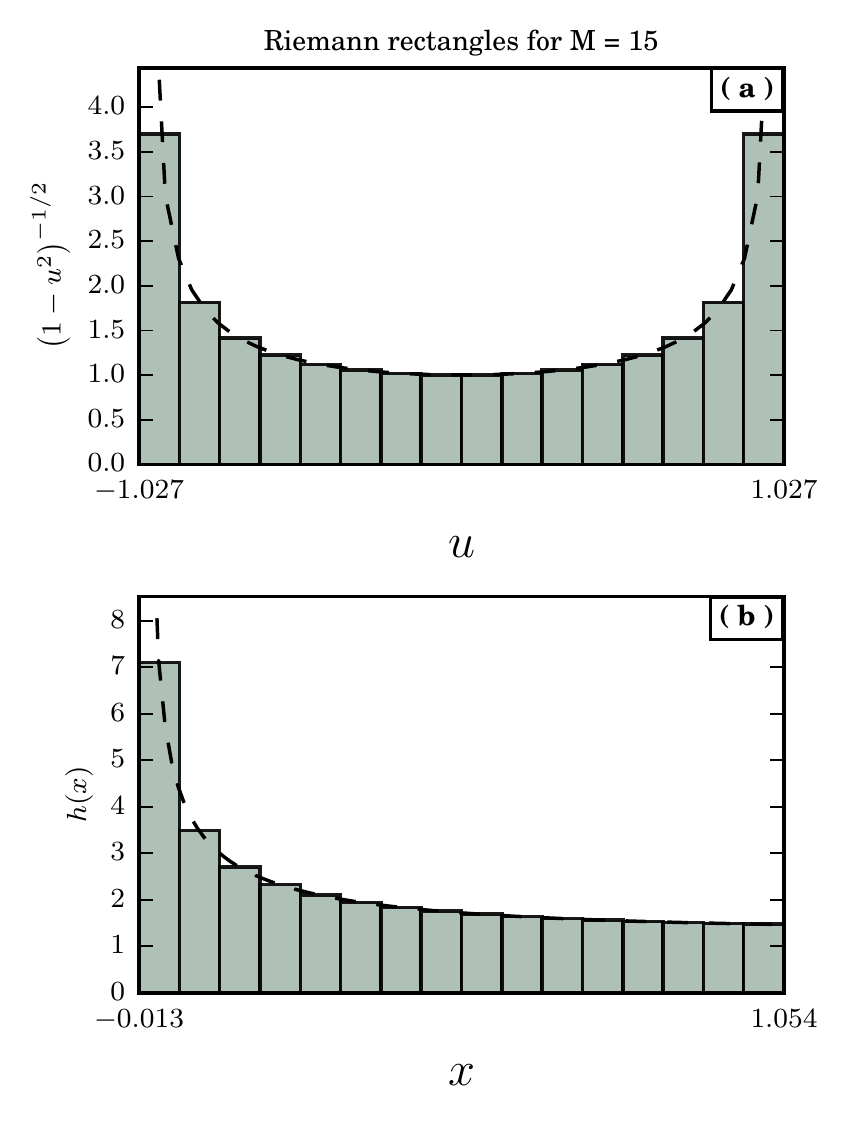}
        \caption{\label{Riemannfig} (Color online) Riemann integrals involved in approximating the sum $B_M$.  (a) The $k$th rectangle has width $\Delta u$ and height $(1-u_k^2)^{-1/2}$. Hence its area represents $A_M(k)$, via \eqref{A_M}. For comparision, the dashed line shows $(1-u^2)^{-1/2}$, which has integrable singularities at $u=\pm 1$. (b) The $k$th rectangle has width $\Delta x = 1/M$ and height $ h(x)=\left(1 - \frac{C_1}{2M} \right) x^{-1/2} + \frac{1}{2}x^{1/2}$. Hence its area approximates $D_{M}^{-}(k)$. } 
\end{figure}

To begin, use the left-right symmetry to combine the $D^{\pm}_M(k)$ sums:
\begin{alignat}{1}
B_M =& \sum_{k=0}^{M}  \left[ A_M(k) - D_{M}^{-}(k) - D_{M}^{+}(k)\right] \notag \\
=& \sum_{k=0}^{M}  \left[ A_M(k) - 2 D_{M}^{-}(k) \right]. \notag 
\end{alignat}
From Fig.~\ref{Riemannfig}(a) and Eq.~\eqref{A_M}, the first sum can be approximated by 
\begin{equation}
\int_{-1}^{+1} (1 - u^2)^{-1/2} du. \notag 
\end{equation}

To handle the second sum, regard 
\begin{equation}
x_k=(C_1/2+k)/M 
\label{x_k}
\end{equation}
for $k=0, \ldots, M$ as mesh points in the midpoint integral for $D^{-}_{M}(k)$, so that $\Delta x = 1/M$. 
Then from Fig.~\ref{Riemannfig}(b) and the first two terms of \eqref{D_M}, we see that the second sum can be approximated by integrating
\begin{equation}
h(x) := \left(1 - \frac{C_1}{2M} \right) x^{-1/2} + \frac{1}{2}x^{1/2}.
\end{equation}
Thus 
\begin{alignat}{1}
B_M = & \int_{-1}^{+1} (1 - u^2)^{-1/2} du. \notag \\ 
&- 2 \int_{0}^{x^*} \left[ \left(1 - \frac{C_1}{2M} \right) x^{-1/2} + \frac{1}{2}x^{1/2}\right] dx\notag \\
&+ O(M^{-2}). \notag 
\end{alignat}

\noindent Here, in the second integral, the upper limit $x^*= x_M+ \Delta x/2= (C_1/2+ M)/M + 1/(2M)$ is the rightmost edge of the rightmost  Riemann rectangle in Fig.~\ref{Riemannfig}(b).  The integral must go up to this value to properly capture the $1/M$ terms.  Evaluating the integrals gives 
\begin{equation}
B_M =\pi - \frac{14}{3} +\frac{ C_1- 3}{2M} + O(M^{-2}). \label{FEq}
\end{equation}

\subsection{Fringe}
Next we need to to estimate the fringe sum $F_M = 2 \sum_{k=0}^{M} D_{M}^{-}(k)$. Substitution of \eqref{D_M} for $D_{M}^{-}(k)$ yields   
\begin{alignat}{1}
F_{M}(k) = & \frac{2}{M^{1/2}} \left[1 - \frac{C_1}{2M} \right] \sum_{k=0}^{M} \left( C_1/2+k \right)^{-1/2} \notag \\
 & + \frac{1}{M^{3/2}} \sum_{k=0}^{M}\left( C_1/2+k \right)^{1/2} \notag \\
& + \frac{1}{M^{3/2}} \left[\frac{C_1 - C_1^2- C_2}{2}\right] \sum_{k=0}^{M} \left( C_1/2+k \right)^{-3/2}.
\label{F_M} 
\end{alignat}
All the sums above are variations on $\sum_{k=0}^{M} (C_1/2 +k)^{-s}$, whose asymptotics for large $M$ follow from results obtained in Refs.~\cite{oliver10, vepstas07}. Those authors showed that, in general,
\begin{alignat}{1}
\sum_{k=0}^{M}\left(k + q\right)^{-s} =& \ \zeta(s, q)- \frac{1}{2} \left(M + 1 + q \right)^{-s} \notag \\
&+ \frac{1}{1 - s}\left( M + 1 + q \right)^{1-s} \notag \\
&+ O(M^{-(2+s)})
\end{alignat}
where $\zeta(s, q)$ is the Hurwitz zeta function \eqref{Hurwitz}. For our purposes, we evaluate the above with $q =  C_1/2$ for $s = -1/2, 1/2$, and $3/2$.  
Inserting these results into \eqref{F_M} gives
\begin{widetext}
\begin{alignat}{1}
F_M =& \frac{14}{3} + 2 \zeta\left(\frac{1}{2}, \frac{C_1}{2} \right) M^{-1/2} + \frac{3 - C_1}{2} M^{-1} \notag \\
&+ \left[ \frac{C_1 - C_2 - C_{1}^2}{2} \zeta\left(\frac{3}{2}, \frac{C_1}{2} \right) - C_1 \zeta\left(\frac{1}{2}, \frac{C_1}{2} \right) + \zeta\left(-\frac{1}{2}, \frac{C_1}{2} \right) \right] M^{-3/2} +O(M^{-2}). \label{DEq}
\end{alignat}
\end{widetext}

\subsection{Locking threshold}

\begin{figure}[!htpb]
        \includegraphics[width=0.9\linewidth]{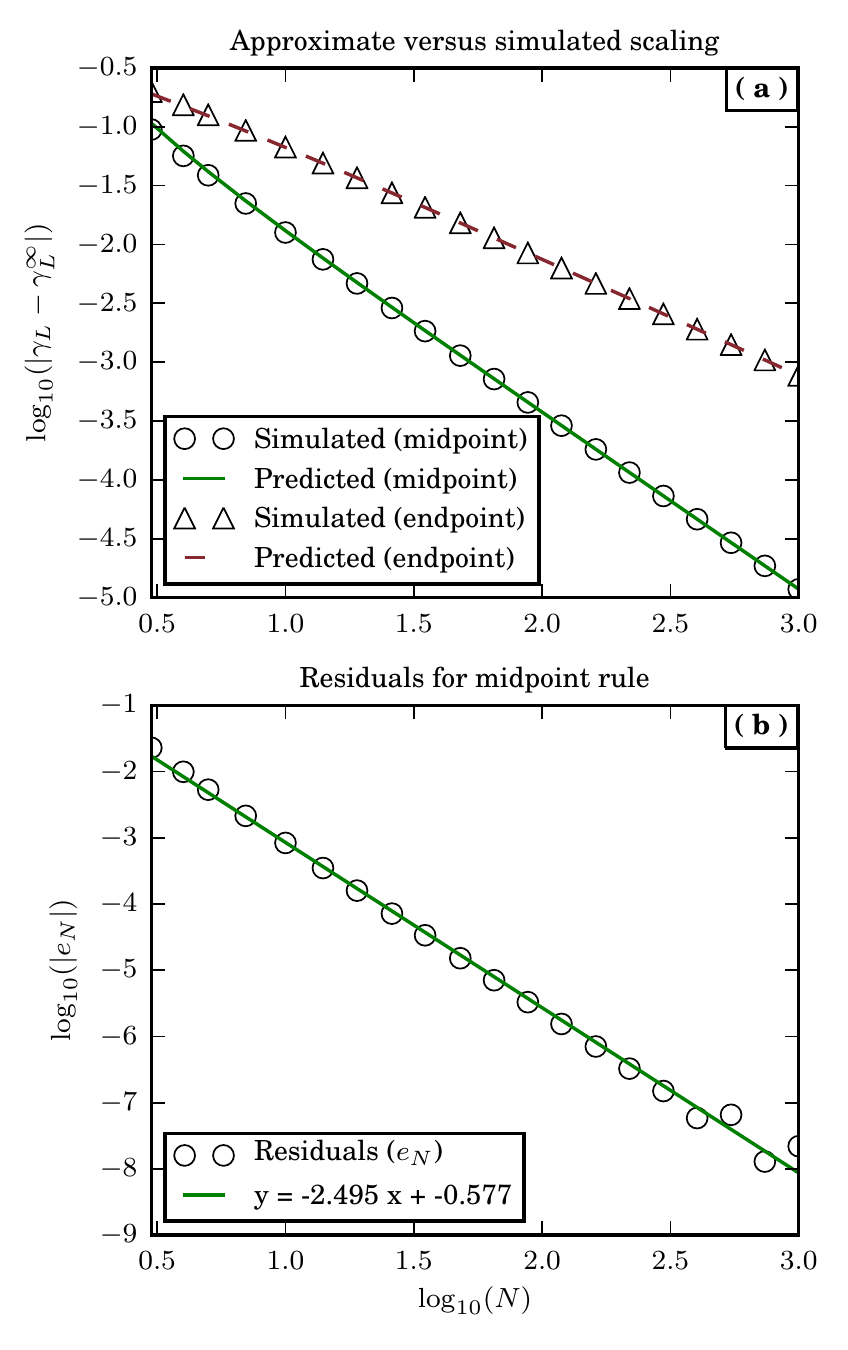}
        \caption{\label{FigCompare} (Color Online) Comparison between theory and numerics. (a) Numerically obtained values of the locking threshold $\gamma_L$ are tested against the predicted scaling laws \eqref{midpoint_prediction} and \eqref{endpoint_prediction} for the midpoint rule \eqref{midpoint} and endpoint rule \eqref{endpoint}, respectively. Numerical values of $\gamma_L$ were obtained by running simulations repeatedly and using a bisection method to locate the bifurcation. Simulations used a fourth-order Runge-Kutta method with a timestep of 0.10, an identification time of $1 \times 10^3$ time units and a maximum transient time of $8 \times 10^4$ time units. Initial phases for each simulation were determined by Eq.~\eqref{sj}.  (b)  Residuals $e_N$, defined as the difference between the predicted and numerical values of $\gamma_L$ for the midpoint rule \eqref{midpoint}. If the prediction \eqref{midpoint_prediction} is correct up to $N^{-3/2}$, then the residual can only be of a higher power, as indeed it appears to be: the measured slope of the least-squares line is close to $-2.5$. The deviations at the right end are numerical artifacts stemming from floating point arithmetic and the finite precision simulations.  }
\end{figure}

With the final piece of the puzzle in place, we can now derive an asymptotic approximation for the locking threshold $\gamma_L$. First recall that for the midpoint rule \eqref{midpoint}, the maximal frequency and the locking threshold are related via $\omega_N= (1-1/N) \gamma_L$, from \eqref{gamma_omega_midpoint}, whereas $\omega_N$ is related to the sum $\alpha_N$ via $\omega_N=(1/4)(1-1/N) \alpha_N$, from \eqref{Riemann_sum_s}; hence 
$\gamma_L= \frac{1}{4}\alpha_N$.
To calculate $\alpha_N = B_M + F_M$, we add \eqref{FEq} and \eqref{DEq}. Then $\gamma_L$ becomes  
\begin{widetext}
\begin{equation}
\gamma_L 
=  \frac{\pi}{4} + \frac{1}{2} \zeta\left(\frac{1}{2}, \frac{C_1}{2} \right) M^{-1/2} \\
 + \frac{1}{4} \left[ \frac{C_1 - C_2 - C_{1}^2}{2} \zeta\left(\frac{3}{2}, \frac{C_1}{2} \right) - C_1 \zeta\left(\frac{1}{2}, \frac{C_1}{2} \right) + \zeta\left(-\frac{1}{2}, \frac{C_1}{2} \right) \right] M^{-3/2} +O(M^{-2}). 
 \notag \\
\end{equation}
\end{widetext}
This equation looks worrisome until one recalls the results of Bailey et al.~\cite{bailey09} mentioned in the Introduction: $C_1$ satisfies $\zeta(1/2, C_1/2) =0$ and $C_2$ is related to $C_1$ by $C_2 = C_1 - C_{1}^2 - 30 \left[\zeta(-1/2, C_1/2)/\zeta(3/2, C_1/2)\right]$.  Then the equation above for $\gamma_L$ collapses to
\begin{equation}
\gamma_L= \frac{\pi}{4}+ 4 \zeta\left(-\frac{1}{2} ,\frac{C_1}{2} \right) N^{-3/2} + O(N^{-2})
\label{midpoint_prediction}
\end{equation}
where we have also made use of $M=N-1$.
The corresponding result for the endpoint rule is 
\begin{equation}
\gamma_L= \frac{\pi}{4} - \frac{\pi}{4}N^{-1} + 4 \zeta\left(-\frac{1}{2} ,\frac{C_1}{2} \right) N^{-3/2} + O(N^{-2}).
\label{endpoint_prediction}
\end{equation}
In these equations, $ 4 \zeta\left(-\frac{1}{2} ,\frac{C_1}{2} \right)\approx 0.3735$, which not only successfully recreates the scaling exponents found numerically \cite{pazo05}, but also gives the exact prefactors. Figure~\ref{FigCompare} shows that these predictions agree with numerics. 

\section{Discussion}

An important thing to notice is that even though the midpoint rule \eqref{midpoint} and the endpoint rule \eqref{endpoint} both approach a uniform distribution on $[-\gamma, \gamma]$, their locking thresholds (their critical widths at which stable phase locking is lost) scale entirely differently. It is possible to blame this on analytic corrections to scaling, but that downplays the level of flexibility we have on the critical exponent. 

For example,  consider the following distribution of natural frequencies:
\begin{equation}
\omega_j = \gamma \left[1 -\frac{\beta}{N^\sigma} \right] \left( 1 - \frac{1}{N} \right) \left[ \frac{2(j-1)}{N-1} \right] \notag
\end{equation}
where $0<\sigma < 3/2$ and $\beta \ne 0$. Then, following our previous analysis, we find to leading order that $\gamma_L = \pi/4 + (\beta\pi/4) N^{-\sigma} + \cdots$. So in principle, the scaling exponent $\sigma$ can be made as close to zero as we like (and the $\beta$ prefactor is equally flexible), despite this distribution also approaching the same uniform distribution as the two we have studied. 

Interestingly, the opposite situation can also happen. Namely, consider the following distribution:
\begin{equation}
\omega_j = \gamma \left[1 - \frac{1}{N} + \frac{16}{\pi} \zeta\left(-\frac{1}{2} ,\frac{C_1}{2} \right) N^{-3/2} \right] \left[ \frac{2(j-1)}{N-1} \right]. \notag
\end{equation}
Numerically speaking, the $(-1/N)$ will dominate the $N^{-3/2}$ term for all $N > 1$, so this still approaches the uniform distribution from the inside. However, if we try to estimate the locking threshold for this case, we get 
\begin{alignat*}{1}
\gamma_L=& \left[1 + \frac{1}{N} - \frac{16}{\pi} \zeta\left(-\frac{1}{2} ,\frac{C_1}{2} \right) N^{-3/2} \right] \left(1 - \frac{1}{N} \right)  \\
& \times \left[\frac{\pi}{4} + 4 \zeta\left(-\frac{1}{2} ,\frac{C_1}{2} \right) N^{-3/2} + O\left(N^{-2} \right) \right] \\
 =& \ \frac{\pi}{4} + O\left(N^{-2}\right).
\end{alignat*}
Now we get a distribution whose locking threshold has a critical exponent larger than what our approach above can estimate. In principle, we could do a more careful analysis using the $C_3$ and $C_4$ coefficients to be able to say exactly what the prefactor and coefficient should be in this case, but we will leave this as an exercise for the diligent student. 

Research supported in part by a Sloan Fellowship to Bertrand Ottino-L\"{o}ffler, and by NSF grants DMS-1513179 and  CCF-1522054.

\small{

}


\begin{thebibliography}{10} 

\bibitem{kuramoto75} Y. Kuramoto, in International Symposium on Mathematical Problems in Theoretical Physics, H. Araki, ed., Lecture Notes in Phys. 39 (Springer, Berlin, 1975), p. 420.

\bibitem{kuramoto84} Y. Kuramoto, Chemical Oscillations, Waves, and Turbulence (Springer, Berlin, 1984).

\bibitem{strogatz00} S. H. Strogatz, Physica D 143, 1 (2000). 

\bibitem{pikovsky03} A. Pikovsky, M. Rosenblum, and J. Kurths, Synchronization: A Universal Concept in Nonlinear Sciences (Cambridge University Press, 2003).

\bibitem{strogatz03} S. H. Strogatz, Sync (Hyperion, New York, 2003).

\bibitem{acebron05} J. A. Acebr\'{o}n, L. L. Bonilla, C. J. P. Vicente, F. Ritort, and R. Spigler, Rev. Mod. Phys. 77, 137 (2005). 

\bibitem{dorfler14} F. D\"{o}rfler and F. Bullo, Automatica 50, 1539 (2014).

\bibitem{pikovsky15} A. Pikovsky and M. Rosenblum, Chaos 25, 097616 (2015)

\bibitem{rodrigues16} F. A. Rodrigues, T. K. DM. Peron, P. Ji and J. Kurths, Phys. Rep. 610, 1 (2016).

\bibitem{winfree67} A. T. Winfree, J. Theor Biol. 16, 15 (1967).


\bibitem{daido87} H. Daido, J. Phys. A 20, L629 (1987).

\bibitem{daido89} H. Daido, Prog. Theor. Phys. 81, 727 (1989)

\bibitem{daido90} H. Daido, J. Stat. Phys. 60, 753 (1990). 

\bibitem{hong05} H. Hong, H. Park, and M. Y. Choi, Phys. Rev. E 72, 036217 (2005). 

\bibitem{hong07} H. Hong, H. Chat\'{e}, H. Park, and L.-H. Tang, Phys. Rev. Lett. 99, 184101 (2007).

\bibitem{hong15} H. Hong, H. Chat\'{e}, L.-H. Tang, and H. Park, Phys. Rev. E 92, 022122 (2015).

\bibitem{hildebrand07} E. J. Hildebrand, M. A. Buice, and C. C. Chow, Phys. Rev. Lett. 98, 054101 (2007).

\bibitem{buice07} M. A. Buice and C. C. Chow, Phys. Rev. E 76, 031118 (2007). 

\bibitem{ermentrout85} G.B. Ermentrout,  J. Math. Biol. 22, 1 (1985).

\bibitem{vanhemmen93} J.L. van Hemmen and W.F. Wreszinski, J. Stat. Phys. 72, 145 (1993). 

\bibitem{jadbabaie04} A. Jadbabaie, N. Motee, and M. Barahona, in Proceedings of the American Control Conference, Boston, MA, 2004, p. 4296. 

\bibitem{maistrenko04} Yu. Maistrenko, O. Popovych, O. Burylko, and P. A. Tass, Phys. Rev. Lett. 93, 084102 (2004).

\bibitem{aeyels04} D. Aeyels and J. Rogge, Prog. Theor. Phys. 112, 921 (2004).

\bibitem{mirollo05} R.E. Mirollo and S.H. Strogatz, Physica D 205, 249 (2005).

\bibitem{pazo05} D. Paz\'{o}, Phys. Rev. E 72, 046211 (2005).

\bibitem{quinn07} D. Quinn, R. Rand, and S. H. Strogatz, Phys. Rev. E 75, 36218 (2007).

\bibitem{verwoerd08} M. Verwoerd and O. Mason, SIAM J. Applied Dynamical Systems 7, 134 (2008). 

\bibitem{verwoerd11} M. Verwoerd and O. Mason, SIAM J. Applied Dynamical Systems 10, 906 (2011).

\bibitem{dorfler11} F. D\"{o}rfler and F. Bullo, SIAM J. Applied Dynamical Systems 10, 1070 (2011). 

\bibitem{bailey09} D. Bailey, J. Borwein, R. Crandall, Experimental Mathematics 18, 107 (2009). 

\bibitem{durgin08} N. Durgin, S. Garcia, T. Flournoy, D. Bailey, \url{https://escholarship.org/uc/item/0kr400n1} (2008).

\bibitem{oliver10} F. W. J. Olver, D. W. Lozier, R. F. Boisvert, and C. W. Clark, NIST Handbook of Mathematical Functions (Cambridge University Press, 2010). 

\bibitem {vepstas07} L. Vep\v{s}tas, \url{http://arxiv.org/abs/math/0702243} (2007). 


\end{thebibliography}
\end{document}